\documentclass[a4paper,11pt]{article}

\usepackage{amsmath,amsfonts,hyperref,amssymb,empheq}
\usepackage{latexsym,mathtools,array}
\usepackage{physics}
\usepackage{bm}
\usepackage[dvipdfmx]{graphicx}

\makeatletter 
\newcommand{\LCASES}[1]{$\m@th\displaystyle{#1}$\hfil}
\newcommand{\CCASES}[1]{\hfil$\m@th\displaystyle{#1}$\hfil}
\newcommand{\RCASES}[1]{\hfil$\m@th\displaystyle{#1}$}
\makeatother
\newcases{ecases}{\quad}{\CCASES{##}}{\LCASES{##}}{\lbrace}{.}
\newcases{ecases*}{\quad}{\CCASES{##}}{{##}\hfil}{\lbrace}{.}

\begin{document}

\title{Log canonical threshold of du Val singularities}
\author{Yoshinori Watanabe \\ \href{mailto::ywatanabe750@fuji.waseda.jp}{ywatanabe750@fuji.waseda.jp}}
\date{\today}
\maketitle

\begin{abstract}
In this paper, we show the log canonical threshold values of the surfaces which has du Val type singularities.These surfaces can be interpreted as statistical or machine learning models.
The results of $A_n, D_n, E_6, E_7$ and $E_8$ are shown.
Bacause the calculation is done on complex fileld, this result would be different from the real fileld, which is in singular learning theory.

\end{abstract} 

\section{Introduction}
  
Most of the models used in statistical and machine learning are singular.
Singular models have degenerated Jacoibians around optimal values.

Their asymtopic behavior are 
descripted by

\section{Log canonical threshold(LCT)}
In \cite{Wata2009}, real log canonical threshold(RCLT) is defined as

\begin{equation}
    \displaystyle
\lambda:=\min_\alpha \min_{1 \leq j \leq d} (\frac{h_j+1}{k_j})  \label{eq:RCLT}
\end{equation}
where $\frac{h_j+1}{k_j}$ is the index of j-th largest pole of zeta function.

\begin{equation}
\displaystyle
\zeta(z)=\int dw K(w)^z\phi(w) =\int_{g^{-1}(C)}K(g(u))^z|g'(u)|du \\
=\sum_\alpha\int_{U_\alpha\cap g^{-1}(C)}u^{kz+h}b(u)du
\end{equation}

$g(u)$ is resolution map of $K(w)$ and $u=\{u_i\}$ are local coordinates of sets $\{U_i\}$.
Then $K(w)$ and $g(u)$ is written as 

$K(g(u))= S u_1^{k_1} u_2^{k_2} \cdots u_d^{k_d}$

$g'(u)= b(u)u_1^{h_1} u_2^{h_2} \cdots u_d^{h_d}$

The index of the pole of zeta function is $\{h_1/k_1,h_2/k_2, \cdots h_d/k_d \}$. for all 

\section{Du val singularities}
Du val singularities is defined on complex filelds. These singularities are enumerated via lie algebric symmetry.
This correspondence is also called or Mckey graph.
And these are depicted by Dynkin diagrams.
\begin{equation*}\begin{split}
A_n: x^2+y^2+z^{n+1}=0 \\
D_n: x^2+y^2z+z^{n-1}=0\\
E_6: x^2+y^3+z^4=0\\
E_7: x^2+y^3+yz^3=0\\
E_8: x^2+y^3+z^5=0\\
\end{split}
\end{equation*}

\section{LCT of Du val singularities}
Here we calculate LCT of each form of Du val singularities.
The following calculation is from \cite{milesreid} and \cite{simplelie}.

\subsection{type $A_n$}
$A_n: x^2+y^2+z^{n+1}=0$
$x,y,z \in \mathbb{C}$

Put 
\begin{equation}
f(x,y,z)=x^2+y^2+z^{n+1} \label{eq:An}
\end{equation}

Then in the region $U_x =\{y_1,z_1\}$ includes $(y,z)=(0,0),x \neq 0$
\begin{equation*}
    g_{x1}=  \begin{ecases*}
        x= & $x_1$ \\
        y= & $y_1x_1$ \\
        z= & $z_1x_1$ 
  \end{ecases*}
\end{equation*}

\begin{equation}
f(x,y,z)=x_1^2(1+y_1^2+x_1^{n-1}z_1^{n+1}) \label{eq:Anx}
\end{equation}

In the region $U_y =\{x_1,z_1\}$ includes $(x,z)=(0,0),y \neq 0$

\begin{equation*}
    g_{y1}=  \begin{ecases*}
        x= & $x_1y_1$ \\
        y= & $y_1$ \\
        z= & $z_1y_1$ 
  \end{ecases*}
  \end{equation*}
  
\begin{equation}
f(x,y,z)=y_1^2(x_1^2+1+y_1^{n-1}z_1^{n+1})  \label{eq:Any}
\end{equation}

\begin{equation*}
    g_{z1}=  \begin{ecases*}
        x= & $x_1z_1$ \\
        y= & $y_1z_1$ \\
        z= & $z_1$ 
\end{ecases*}
\end{equation*}

and in the region $U_z =\{x_1,y_1\}$ includes $(x,y)=(0,0),z \neq 0$

\begin{equation}
f(x,y,z)=z_1^2(x_1^2+y_1^2+z_1^{n-1}) \label{eq:Anz}
\end{equation}

respectively. The formulars \eqref{eq:Anx}, \eqref{eq:Any} for $U_x$ and $U_y$ are not singlar.
The second factor of \eqref{eq:Anz} is same as the form of $A_{n-2}$

\begin{equation}
x_1^2+y_1^2+z_1^{n-1} \label{eq:Anm2}
\end{equation}

The resolution process applied recursively . $A_n$ can be written as
\begin{equation}
    z_k^{n+1} (x_1^2+y_1^2+1) \label{eq:Anneven}
\end{equation}
(even case, k=(n+1)/2) or 
\begin{equation}
    z_k^{n} (x_1^2+y_1^2+z) \label{eq:Annodd}
\end{equation}
(odd case, k=n/2). By appling strict transformations$\{g_{z_i}\}_{i\leq k}$
\begin{equation*}
    g_{z_i}= \begin{ecases*}
        x= & $x_1z_{i-1}=x_1z_i$ \\
        y= & $y_1z_{i-1}=y_1z_i$ \\
        z_{i-1}= & $z_i$ 
\end{ecases*}
\end{equation*}
  
$g:=g_{z_n} \circ g_{z_{n-1}} \circ \cdots \circ g_{z_1}$

The determinant Jacoibian of g is 

\begin{equation}
\begin{split}
|g'|&:=|\partial g(z_n)/ \partial z_1|\\
&=|\partial (g_{z_k} \circ g_{z_{k-1}} \circ \cdots \circ g_{z_1})/ \partial z_1|\\
&=|(\partial g_{z_k} / \partial z_k)||(\partial g_{z_{k-1}} / \partial z_{k-1})| \cdots |(\partial g_{z_2} / \partial z_2)|(\partial g_{z_{1}} / \partial z_{1})|\\
&=z_k^2 z_{k-1}^2 \cdots z_{2}^2  z_{1}^2 |\\
&=z_k^{2k}
\end{split}
\end{equation}

The LCT of \eqref{eq:Anneven} is $\lambda=\frac{2k+1}{n+1}=\frac{n+2}{n+1}$,

LCT of \eqref{eq:Annodd} is  $\lambda=\frac{2k+1}{n}=\frac{n+1}{n}$.

\subsection{type $D_4$}
\begin{equation}
f(x,y,z)=x^2+y^2z+z^{3}
\end{equation}

$x^2+z(y-iz)(y+iz)=:x^2+abc$
$(a,b,c \in \mathbb{C})$

This can be converted to $a',b',c' \in \mathbb{C}$ and new $x,y,z$ which

$x^2+a'b'c'=x^2(y+z)(y+jz)(y+j^2z)$ , $(j=e^{2\pi i/3})$

by linear conversion. 

Then strict transformation
\begin{equation*}
    g_x= \begin{ecases*}
        x= & $x_1$\\
        y= & $x_1y_1 =xy_1$\\
        z= & $x_1z_1 =xz_1$\\
\end{ecases*}
\end{equation*}
 makes 
 $x^2(1+x(y_1^3+z_1^3))$. This is not singular
 
 Next 
\begin{equation*}
g_y= \begin{ecases*}
     x= & $x_1y_1 =x_1y$\\
     y= & $y_1$\\
     z= & $y_1z_1 =z_1y$\\
\end{ecases*}
\end{equation*}
is applied and obtain $y^2(x_1+y(1+z_1^3))$. 
The exceptional curve is $E:y = x_1=0$ and it intersects three $A_1$ singularities at $(0,0,a)$ with $a^3=-1$.

For z ($g_z$)is same as the case of y.
As a result LCT $\lambda=\frac{2+1+1}{2+1}=\frac{4}{3}$

\subsection{type $D_5$}

$f(x,y,z)=x^2+y^2z+z^{4}=0$

For variable $y$ strict transformation
\begin{equation*}
    g_y= \begin{ecases*}
         x= & $x_1y_1 =x_1y$\\
         y= & $y_1$\\
         z= & $y_1z_1 =z_1y$\\
\end{ecases*}
\end{equation*}

$y^2(x_1^2+yz_1+y^2z_1^4)$.
This is converted to $y^2(x_1^2+y_1z_1')$ by variable transformations $z_1'=z_1+yz_1^4,y=y,x=x$ 
whose Jacoibian determinant is $|g'|=1$.
The formular inside brackets is same as $A_1: x^2+y^2+z^2$. So the LCT at this local chart is $\frac{3+2}{2+2}$.

On the other hand, with strict transformation
\begin{equation*}
    g_z= \begin{ecases*}
         x= & $x_1z_1 =x_1z$\\
         y= & $y_1z_1 =y_1z$\\
         z= & $z_1$\\
\end{ecases*}
\end{equation*}
$z^2(x_3^2+z+zy_3^3+z^2)$.
This can convert $f(x,y,z)$ to    
$z^2(x_3^2+(z+\frac{y^2}{2})^2-\frac{y^4}{4})$
by lienar conversion. And this is same form as $A_3$, whose $\lambda=\frac{4}{3}$.

As a result the LCT of $D_5$ is $\lambda=\min\{\frac{5}{4},\frac{6}{5}\}=\frac{6}{5}$

\subsection{type $D_n$}
\begin{equation}
D_n: x^2+y^2z+z^{n-1}=0
\end{equation}

For the chart $U_x$, strict transformation $g_x$ makes
$x^2(1+xy_1^2z_1+x^{n-3}z_1^{n-1})$
all factors are normal crossing and this is not singular.

For the chart $U_y$, strict transformation $g_y$ makes $y^2(x_2^2+yz_2+y_2^{n-3}z_2^{n-1})$.
One more blow-up is enough for resolving this singularity.
For example, $g_{z_3}$ makes $y_2^2z_2^2(x_3^2+y_3+y_3^{n-3}z_3^{2n-5})$
The gradient is $\nabla (2x_3,1+(n-3)y_3^{n-4}z_{2n-5}, (2n-5)y_3^{n-3}z_3^{2n-6}) \neq 0$
 at (0,0,0).
 
For the chart $U_z$, strict transformation $g_z$ makes 
$z_1^2(x_1^2+y_1^2z_1+z_1^{n-2})$. This is the form of $D_{n-2}$ and resolution process is applied recursively.
Then total factor for 
By using the result of $D_4$ and $D_5$, total LCT value is $\lambda=\frac{n+1}{n} or frac{n-1}{n-2}$.

\subsection{type $E_6$}
\begin{equation}
E_6: x^2+y^3+z^4=0
\end{equation}
resolution for $U_y$ by $g_y$ converts $E_6$ to $y_1^2(x_1^2+y_1^2+y_1^2z_1^4)$. This can be conert to $y_1^2(x_1^2+y_1^2z_2^4)$ by linear transformation and become 
\begin{equation}
Y(x^2+YZ) \label{eq:E62y}
\end{equation}
by $Y=y^2,Z=z^4$ which is same as $A_1$.
In case $U_z$, the formular is converted to 
\begin{equation}
z_1^2(x_1^2+z_1y_1^3+z_1^2)
\end{equation}
As second step, For x and z  this is smooth curve $z_2^2x_2^2(1+x_2^2z_2y_2^3+z_2^3)$ and for y 
\begin{equation}
z_2^2y_2^4(x_2^2+y_2^2z_2+z_2^2)
\end{equation}
As third step, 
\begin{equation}
    z_3^2y_3^{6+2}(x_3^2+y_3z_3+z_3^2) \label{eq:E63y}
\end{equation}
\begin{equation}
    x_3^2z_3^2y_3^4x_3^4x_3^2(1+x_3y_3-2z_3+z_3^2)\label{eq:E63z}
\end{equation}
As forth step, \eqref{eq:E63y} converted to 
\begin{equation}
    y_4^2z_4^2y_4^8y_4^2(x_4+z_4+z_4^2)\label{eq:E64y}
\end{equation}
or 
\begin{equation}
    x_4^{10}z_3^2y_4^8x_4^2(1+y_4z_4+z_4^2)\label{eq:E64z}
\end{equation}

The LCT is minumum for all of \eqref{eq:E62y}, \eqref{eq:E63z}, \eqref{eq:E64y} and \eqref{eq:E64z}.
\begin{equation}
\begin{split}    
\lambda&=\min\{ \frac{1+2+2}{1+2} , 
\{ \infty, \frac{12+1}{12}, \frac{2+1}{2} \},
\{ \frac{12+1}{12}, \frac{8+1}{8}, \frac{2+1}{2} \},
\{ \frac{2+1}{8}, \frac{6+1}{4}, \frac{2+1}{2} \}
\}\\
&=\frac{12}{13}\\
\end{split}
\end{equation}

\subsection{type $E_7$}
\begin{equation}
E_7: x^2+y^3+yz^3=0
\end{equation}
For $U_z$, applying
$g_z$
\begin{equation*}
    g_z= \begin{ecases*}
         x= & $x_1z_1 =x_1z$\\
         y= & $y_1z_1 =y_1z$\\
         z= & $z_1$\\
\end{ecases*}
\end{equation*}

$z_1^2(x_1^2+z_1y_1^3+)$
The formular inside brackets is the form of $D_6$.
Because $zy^3+yz^3=yz(z+y^2)\propto w(v-iw^2)(v-iw^2+2iw^2)=w(v-iw^2)(v+iw^2)-w(v^2+w^4)$ by putting $z=(v-iw), y^2=2iw^2$.
$x_1^2 wv^2+w^5$ .
The Jacobian of convertion $(x,y,z)\rightarrow (x,v,w)$ is constant. Then $z_1^2w^2(x^2+wv^2+w^3)$ is converted to 
$(v^2+w^2)w^2x^2(1+wxv+xw^3)$  or $(v^2+w^2)w^4(x^2+wxv+xw^3)$ for local coordinate $U_v$ or $U_w$.

Total LCT is $\lambda=\min\{ \frac{2+1}{2},\frac{2+2+1}{4},\frac{2+2+1}{6}\}=\min\{\frac{3}{2},\frac{5}{4},\frac{5}{6}\}=\frac{5}{6}$ which is minimum of the ratio of index of x,v and w.

\subsection{type $E_8$}
\begin{equation}
E_8: x^2+y^3+z^5=0
\end{equation}
Appliyng $g_{z_1}$,the formular becomes $z_1^2(x_1^2+zy_1^3+y_1z_1^3)$. Because this is the form as $E_7$, the RCT is $\lambda=\frac{1+2+6}{2+6}=\frac{9}{8}$

\section{Discussion and Limitation}
The result suggests that LCT of degree n singular curve on complex field is $\frac{n+1}{n}$.This is different from the one of real, not algebric closed field case.
RCLT is coefficient of loss fuction of machine learning, or a posterior distribution of statistical model.
And much complicated treatment  even two or three layer neural network as shown in  \cite{Farrugia-Roberts2022}, \cite{wong2022}.
LCT on complex filed may have some meaning in quantum computing, quantum information.

\bibliographystyle{jplain}
\bibliography{Duval_ref}
\end{document}